
\font\bigbf=cmbx10 scaled \magstep3
\font\medbf=cmbx10 scaled \magstep2

\magnification=\magstep1
\null
\vskip 3truecm
\centerline{\bigbf Extremal Selections of Multifunctions}
\vskip 4truemm
\centerline{\bigbf Generating a Continuous Flow}
\vskip 15truemm
\centerline{\sl Alberto Bressan and Graziano Crasta}
\vskip 2em
\centerline{S.I.S.S.A. -- Via Beirut 4, Trieste 34014, ITALY}
\vskip 2cm

\magnification=\magstep1
\def \n{\noindent}
\def \v{\vskip 1em}

\def \vsk{\vskip 4em}
\def \rn {I\!\!R^n}

\def \M{{\cal M}}
\def \Me{{\cal M}^{ext}}
\def \yij{y^i_j}
\def \O{\Omega^\dagger}
\def \sus {(s,u(s))}
\def\dom {[0,T]\times\Omega}
\def \R {I\!\!R}
\def \txi{(t_i,x_i)}
\def \fiij{\varphi_j^i}

\def \sumj{\sum_{j=0}^n}
\def\fjtx{\varphi_j^{(t,x)}}
\def\pjtx{\psi_j^{(t,x)}}
\def\psiij{\psi_j^i}
\def \AC{\cal AC}
\def \S{{\cal S}}
\def \ve{\varepsilon}
\def \P {{\cal P}}

\def \inda {(t_0,x_0)\in [0,T]\times D}
\def \L {{{\cal L}^1}}
\def \i{\item}
\def \q{\rightline {Q.E.D.}}
\baselineskip = 16 pt
\hfuzz=24pt
\n
{\medbf 1 - Introduction}
\v
Let $F:[0,T]\times\R^n\mapsto 2^{\R^n}$ be a continuous multifunction with
compact, not necessarily convex values.  If $F$ is Lipschitz continuous,
it was shown in [4] that
there exists a measurable selection $f$ of $F$ such that, for every
$x_0$, the Cauchy problem
$$\dot x(t)=f(t,x(t)),\qquad\qquad x(0)=x_0$$
has a unique Caratheodory solution, depending continuously on $x_0$.

In this paper, we prove that the above selection $f$ can be chosen so that
$f(t,x)\in extF(t,x)$ for all $t,x$.  More generally, the result remains
valid if $F$ satisfies the following Lipschitz Selection Property:
\v
\i{(LSP)} {\sl For every $t,x$, every $y\in \overline{co} F(t,x)$ and
$\ve>0$, there exists a Lipschitz selection $\phi$ of $\overline{co}F$,
defined on a neighborhood of $(t,x)$, with $|\phi(t,x)-y|<\ve$.}
\v
We remark that, by [7,9], every Lipschitz multifunction with compact
values satisfies (LSP).  Another interesting class, for which
(LSP) holds, consists of those continuous multifunctions $F$ whose values are
compact and have convex closure with nonempty interior.  Indeed, for any given
$t,x,y,\ve$, choosing $y'\in int\, \overline{co}F(t,x)$ with $|y'-y|<\ve$,
the constant function $\phi\equiv y'$ is a local selection from
$\overline{co}F$ satisfying the requirements.

In the following, $\Omega\subseteq\R^n$ is an open set, $\overline B(0,M)$
is the closed ball centered at the origin with radius $M$,~ $\overline
B(D;~MT)$ is the closed neighborhood of radius $MT$ around the set $D$,
while $\AC$ the Sobolev space of all
absolutely continuous functions $u:[0,T]\mapsto\R^n$, with norm
$\|u\|_{\AC}=\int_0^T\big(|u(t)|+|\dot u(t)|\big)~dt$.
\v
\n{\bf Theorem 1.}  {\it Let $F:[0,T]\times \Omega\mapsto
2^{\R^n}$ be a bounded continuous multifunction with compact values,
satisfying (LSP).  Assume that $F(t,x)\subseteq\overline B(0,M)$ for all $t
,x$ and let $D$ be a compact set such that $\overline B(D;~MT)\subset
\Omega$.
Then there exists a measurable function
$f$, with
$$f(t,x)\in ext F(t,x)\qquad\qquad\forall t,x,\eqno(1.1)$$
such that, for every $(t_0,x_0)\in [0,T]\times D$,
the Cauchy problem
$$\dot x(t)=f(t,x(t)),\qquad\qquad x(t_0)=x_0\eqno(1.2)$$
has a unique Caratheodory solution $x(\cdot)=x(\cdot,t_0,x_0)$ on $[0,T]$,
depending continuously on $t_0,x_0$ in the norm of $\AC$.

Moreover, if $\ve_0>0$ and a Lipschitz continuous selection $f_0$
of $\overline{co}F$ are given, then one can construct $f$ with the
following additional
property.  Denoting by $y(\cdot, t_0,x_0)$ the unique solution of
$$\dot y(t)=f_0(t,y(t)),\qquad\qquad y(t_0)=x_0,\eqno(1.3)$$
for every $(t_0,x_0)\in [0,T]\times D$ one has
$$\big| y(t,t_0,x_0)-x(t,t_0,x_0)\big|\leq\ve_0\qquad\qquad\forall t
\in [0,T].\eqno(1.4)$$}
\v
The proof of the above theorem, given in section 3, starts
with the construction of a
sequence $f_n$ of selections from $\overline{co}F$, which are piecewise
Lipschitz
continuous in the $(t,x)$-space.  For every $u:[0,T]\mapsto \R^n$ in
a class of Lipschitz continuous functions, we then show that the composed maps
$t\mapsto f_n(t,u(t))$ form a Cauchy sequence in $\L\big([0,T];~\R^n\big)$,
converging pointwise almost everywhere to a map of the form $f(\cdot, u
(\cdot))$, taking values within the extreme points of $F$.
This convergence is obtained through an argument which is considerably
different from previous works. Indeed, it relies on a careful use of
the likelihood functional introduced in [3], interpreted here
as a measure of ``oscillatory non-convergence" of a set of derivatives.

Among various corollaries, Theorem 1 yields an extension, valid for the
wider class of multifunctions with the property (LSP),
of the following results, proved in [5], [4] and [6], respectively.
\v
\i{(i)} Existence of selections from the solution set of a differential
inclusion, depending continuously on the initial data.
\v
\i{(ii)} Existence of selections from a multifunction, which generate a
continuous flow.
\v
\i{(iii)} Contractibility of the solution sets of $\dot x\in F(t,x)$
and $\dot x\in ext F(t,x)$.
\v
These consequences, together with an application to
bang-bang feedback controls, are described in section 4.
\vsk

\n{\medbf 2 - Preliminaries}
\v
As customary, $\bar{A}$\ and $\overline{\rm co}\, A$
{}~denote here the closure and the closed convex hull of $A$ respectively,
while $A\backslash B$\ indicates a set--theoretic difference.
The Lebesgue measure of a set $J\subset\R$\ is $m(J)$.  The
characteristic function of a set $A$\ is written as $\chi_{\strut A}$.

In the following, ${\cal K}_n$\ denotes the family of all
nonempty compact convex subsets of $\rn$, endowed with Hausdorff
metric. A key technical tool used in our proofs will be the
function $h:\R^n\times {\cal K}_n\mapsto \R\cup\{-\infty\}$,
defined by
$$h(y,K)\doteq\sup\left\{
\left(\int_0^1 |w(\xi)-y|^2\,d\xi\right)\!^{1\over 2};\quad
w:[0,1]\mapsto K,~~~\int_0^1 w(\xi)\,d\xi =y\right\}\eqno(2.1)$$
with the understanding that $h(y,K)=-\infty$\ if $y\not\in K$.
Observe that $h^2(y,K)$\ can be interpreted as the maximum variance
among all random variables supported inside $K$, whose mean
value is $y$. The following results were proved in [3]:
\v
\n{\bf Lemma 1.}  {\it
The map $(y,K)\mapsto h(y,K)$\ is upper semicontinuous in both variables;
for each fixed $K\in{\cal K}_n$ the function
$y\mapsto h(y,K)$\ is strictly concave down on $K$. Moreover, one has
$$h(y,K)=0\qquad\hbox{ if and only if}\qquad y\in ext K,\eqno(2.2)$$
$$h^2(y,K)\leq r^2(K)-\big|y-c(K)\big|^2,\eqno(2.3)$$
where $ c(K) $ and $r(K)$ denote the Chebyschev center and the Chebyschev
radius of $K$, respectively.}
\v
For the basic theory of multifunctions and differential inclusions
we refer to [1].  As in [2], given a map $g:[0,T]\times \Omega\mapsto\rn$,
we say that $g$\ is directionally continuous along the directions of the
cone $\Gamma^N=\big\{(s,y); ~|y|\leq Ns\big\}$~ if
$$g(t,x)=\lim_{k\rightarrow\infty}g(t_k,x_k)$$
for every $(t,x)$ and every sequence $(t_k,x_k)$\ in the domain of $g$\
such that
$t_k\rightarrow t$\ and $|x_k-x|\leq N(t_k-t)$~
for every $k$.  Equivalently, $g$ is $\Gamma^N$-continuous iff it is
continuous w.r.t. the topology
generated by the family of all conical neighborhoods
$$\Gamma^N_{(\hat t,\hat x,\varepsilon)}\doteq
\big\{(s,y)\ ;~~~\hat t\leq s\leq \hat t+\varepsilon,\ |y-\hat x|\leq N(s-t)
\big\}.\eqno(2.4)$$
A set of the form (2.4) will be called an $N$--cone.

Under the assumptions on $\Omega,~D$ made in Theorem 1, consider the
set of Lipschitzean functions
$$Y\doteq\big\{ u:[0,T]\mapsto\overline B(D,MT);\qquad |u(t)-u(s)|
\leq M|t-s|\quad\forall t,s\big\}.$$  The
Picard operator of a map $g:[0,T]\times\Omega\mapsto\R^n$ is defined as
$${\cal P}^g(u)(t)\doteq \int_0^t g\sus ~ds\qquad\qquad u\in Y.$$
The distance between two Picard operators will be measured by
$$\big\|\P^f-{\cal P}^g\big\|\ =\ \sup \left\{
\left|\int_0^t [f\sus - g\sus]\;ds\right|\ ;~~~\ t\in [0,T],~~ u\in Y
\right\}.\eqno(2.5)$$
The next Lemma will be useful in order to prove the uniqueness
of solutions of the Cauchy problems (1.2).
\v
\n {\bf Lemma 2.}  {\it Let $f$ be a measurable map from
$\dom$ into $\overline B(0,M)$,
with $\P^f$
continuous on $Y$. Let $D$ be compact, with $\overline B(D,~
MT)\subset\Omega$, and assume that the Cauchy problem
$$\dot x(t)=f(t,x(t)),\qquad x(t_0)=x_0,\qquad~~~t\in[0,T]\eqno(2.6)$$
has a unique solution, for each $\inda$.

Then, for every $\epsilon>0$,\ there exists $\delta>0$\
with the following property.
If $g:\dom\rightarrow\overline B(0,M)$\ satisfies
$\big\|\P^g-\P^f\big\|\leq\delta$, then for every $(t_0,x_0)\in [0,T]\times
D$, any solution of the Cauchy problem
$$\dot y(t)=g(t,y(t))\qquad y(t_0)=x_0\qquad~~~t\in [0,T]\eqno(2.7)$$
has distance $<\ve$ from the corresponding solution of (2.6).
In particular, the solution set of (2.7) has diameter $\leq 2\ve$
in ${\cal C}^0\big([0,T];\R^n\big)$.}
\v
{\bf Proof.}  If the conclusion fails, then there exist sequences
of times $t_\nu,~t'_\nu$,
maps $g_\nu$ with $\big\| \P^{g_\nu}-\P^f\big\|\to 0$, and
couples of solutions $x_\nu,
y_\nu:[0,T]\mapsto \overline B(D;MT)$ of
$$\dot x_\nu(t)=f(t,x_\nu(t)),\qquad\quad
\dot y_\nu(t)=g_\nu(t,y_\nu(t))\qquad\qquad t\in [0,T],\eqno(2.8)$$
with $$x_\nu(t_\nu)=y_\nu(t_\nu)\in D,\qquad\qquad
\big|x_\nu(t'_\nu)-y_\nu(t'_\nu)\big|\geq\ve\qquad\quad
\forall\nu.\eqno(2.9)$$
By taking subsequences, we can assume that $t_\nu\to t_0$,
$t'_\nu\to \tau$, $x_\nu(t_0)\to x_0$,
while $x_\nu\to x$ and  $y_\nu\to y$ uniformly
on $[0,T]$.  From (2.8) it follows
$$\eqalign{&\left| y(t)-x_0-\int_{t_0}^t f(s,y(s))~ds\right|
\leq\big|y(t)-y_\nu(t)\big|+\big|x_0-y_\nu(t_0)\big|\cr
&\qquad\qquad +\left|\int_{t_0}^t\big[ f(s,y(s))-f(s,y_\nu(s))\big]
{}~ds\right|+\left|\int_{t_0}^t\big[ f(s,y_\nu(s))-g_\nu(s,y_\nu(s))
\big]~ds\right|.\cr}\eqno(2.10)$$
As $\nu\to\infty$, the right hand side of (2.10) tends to zero,
showing that $y(\cdot)$ is a solution of (2.6).
By the continuity of $\P^f$,~ $x(\cdot)$ is also a solution of
(2.6), distinct from $y(\cdot)$ because
$$|x(\tau)-y(\tau)|=\lim_{\nu\to\infty} |x_\nu(\tau)-y_\nu(\tau)|
=\lim_{\nu\to\infty}\big|x_\nu(t'_\nu)-y_\nu(t'_\nu)\big|
\geq\ve.$$
This contradicts the uniqueness assumption, proving the lemma.
\vsk
\n{\medbf 3 - Proof of the main theorem}
\v
Observing that $ext F(t,x)=ext \overline{co}F(t,x)$ for every compact set
$F(t,x)$, it is clearly not restrictive to prove Theorem 1 under the
additional assumption that all values of $F$ are convex.  Moreover,
the bounds on $F$ and $D$ imply that no solution of the Cauchy problem
$$\dot x(t)\in F(t,x(t)),\qquad\quad x(t_0)=x_0,\qquad\quad t\in [0,T],$$
with $x_0\in D$, can escape from the set $\overline B(D,~MT)$.  Therefore,
it suffices to construct the selection $f$ on the compact set
$\O\doteq [0,T]\times\overline B(D,~MT)$.  Finally, since every
convex valued multifunction satisfying (LSP) admits a globally defined
Lipschitz selection, it suffices to prove the second part of the theorem,
with $f_0$ and $\ve_0>0$ assigned.

We shall define a sequence of directionally continuous selections of $F$,
converging a.e. to a selection from $extF$.
The basic step
of our constructive procedure will be provided by the next lemma.
\v
\n{\bf Lemma 3.}  {\it Fix any $\ve>0$. Let $S$ be a compact subset
of $[0,T]\times\Omega$ and let
$\phi:S\rightarrow\rn$\ be a continuous selection of $F$ such that
$$h\big(\phi(t,x),F(t,x)\big)<\eta\qquad\qquad \forall (t,x)\in S,
\eqno(3.1)$$ with $h$ as in (2.1).
Then there exists a piecewise Lipschitz selection
$g:S\rightarrow\rn$\ of $F$\ with the following properties:
\v
\item{(i)} There exists a finite covering $\{\Gamma_i\}_{i=1\ldots,\nu}$,
consisting of $\Gamma^{M+1}$--cones, such that, if we define the
pairwise disjoint sets
$\Delta^i \doteq \Gamma_i\setminus \bigcup_{\ell<i} \Gamma_\ell$, then
on each $\Delta^i$\ the following holds:
\v
\itemitem{(a)} there exist Lipschitzean selections $\psi^i_j :
\overline{\Delta^i}\mapsto\rn$, ~$j=0,\ldots,n$, such that
$$g\big|_{\Delta^i} = \sum_{j=0}^n \psi^i_j\ \chi_{\strut A^i_j}.
\eqno(3.2)$$
where each $A^i_j$ is a finite union of strips of the form
$\big([t',t'')\times \rn\big)\cap\Delta^i$.
\v
\itemitem{(b)}  For every $j=0,\ldots, n$\ there exists an affine map
$\fiij(\cdot)=\langle a^i_j,\cdot\rangle +b^i_j$ such that
$$\fiij\big(\psi^i_j(t,x)\big)\leq\ve,\qquad
\fiij(z)\geq h(z,F(t,x)),\qquad\qquad\forall (t,x)\in\overline{\Delta^i}
,~~z\in F(t,x).\eqno(3.3)$$
\v
\item{(ii)}  For every $u\in Y$ and every interval $[\tau,\tau']$
such that $(s,u(s))\in S$ for $\tau\leq s<\tau'$, the
following estimates hold:}
$$\left| \int_\tau^{\tau'} \big[\phi\sus - g\sus\big]~ds\right |\ \leq \ve,
\eqno(3.4)$$
$$\int_\tau^{\tau'}\big| \phi\sus-g\sus\big|~ds\ \leq\ \ve+\eta(\tau'-\tau).
\eqno(3.5)$$
\v
\n{\bf Remark 1.}  Thinking of $h(y,K)$ as a measure for the distance of $y$
from the extreme points of $K$, the above lemma can be interpreted
as follows.  Given any selection $\phi$ of $F$, one can find a $\Gamma^{M+
1}$-continuous selection $g$ whose values lie close to
the extreme points of $F$ and whose Picard operator $\P^g$, by (3.4),
is close to $\P^\phi$.   Moreover, if the values of
$\phi$ are near the extreme points of $F$, i.e. if $\eta$ in (3.1) is
small, then $g$ can be chosen close to $\phi$.  The estimate (3.5)
will be a direct consequence of the definition (2.1) of $h$ and of
H\"older's inequality.
\v
\n{\bf Remark 2.}  Since $h$ is only upper semicontinuous, the
two assumptions
$y_\nu\to y$ and $h(y_\nu,K)\to 0$ do not necessarily imply $h(y,K)=0$.
As a consequence, the a.e. limit of a convergent sequence of approximately
extremal selections $f_\nu$ of $F$ need not take values inside $ext F$.
To overcome this difficulty, the estimates in (3.3) provide
upper bounds for
$h$ in terms of the affine maps $\fiij$.  Since each $\fiij$ is continuous,
limits of the form $\fiij (y_\nu)\to\fiij (y)$ will be straightforward.
\v
\n{\bf Proof of Lemma 3.\ } For every $(t,x)\in S$\ there exist values
$y_j(t,x)\in F(t,x)$ and coefficients $\theta_j(t,x)\geq 0$, with
$$\phi(t,x)=\sumj \theta_j(t,x)y_j(t,x),\qquad\qquad\sumj
\theta_j(t,x)=1,$$
$$h\big(y_j(t,x),F(t,x)\big)<\ve/2.$$
By the concavity and the upper semicontinuity of $h$,
for every $j=0,\ldots,n$\ there exists an affine function
$\fjtx(\cdot)=\langle a_j^{(t,x)},\cdot\rangle +b_j^{(t,x)}$
such that
$$\fjtx(y_j(t,x))<h\big(y_j(t,x),F(t,x)\big)+{\ve\over 2}<\ve,$$
$$\fjtx(z)>h\big(z,F(t,x)\big)\qquad\qquad\forall z\in F(t,x).$$
By (LSP) and the continuity of each $\fjtx$,
there exists a neighborhood ${\cal U}$\
of $(t,x)$\ together with Lipschitzean selections
$\pjtx:{\cal U}\mapsto\rn$,
such that, for every $j$ and every $(s,y)\in {\cal U}$,
$$\Big|\pjtx(s,y)-y_j(t,x)\Big|<{\ve\over
4T},\eqno(3.6)$$
$$\varphi_j^{(t,x)}\big(\pjtx(s,y)\big)<\ve.\eqno(3.7)$$
Using again the upper semicontinuity of $h$, we can find a neighborhood
${\cal U}'$\ of $(t,x)$\ such that
$$\fjtx(z)\geq h\big(z,F(s,y)\big)\quad
\qquad\forall z\in F(s,y),~~~(s,y)\in {\cal U}',~~~j=0,\ldots,n.\eqno(3.8)$$
Choose a neighborhood $\Gamma_{t,x}$ of $(t,x)$, contained in
${\cal U}\cap{\cal U}'$, such that, for every
point $(s,y)$ in the closure $\overline \Gamma_{t,x}$, one has
$$|\phi(s,y)-\phi(t,x)|<{\ve\over 4T},\eqno(3.9)$$
It is not restrictive to assume that $\Gamma_{t,x}$ is a
$(M+1)$-cone, i.e. it has the form (2.4) with $N=M+1$.
By the compactness of $S$\ we can extract a finite subcovering
$\big\{\Gamma^i;~~1\leq i\leq\nu\big\}$, with $\Gamma_i\doteq\Gamma_{t_i,x_i}$.
Define $\Delta^i\doteq\Gamma_i\setminus\bigcup_{j<i}\Gamma_j$\ and
set $\theta_j^i=\theta_j\txi$, ~$y_j^i=y_j\txi$, ~$\psiij=\psi_j^{\txi}$,
{}~$\fiij={\varphi_j}^{\txi}$.
Choose an integer $N$\ such that
$$N>{8M\nu^2T\over\ve}\eqno(3.10)$$
and divide $[0,T]$\ into $N$\ equal subintervals $J_1,\ldots,J_N$,
with
$$J_k=\big[t_{k-1},~t_k\big),\qquad\qquad t_k={kT\over N}.\eqno(3.11)$$
For each $i,k$ such that $\big(J_k\times\rn\big)\cap\Delta^i
\not= \emptyset$, we then split $J_k$ into $n+1$ subintervals
$J^i_{k,0},\ldots, J^i_{k,n}$ with lengths proportional to
$\theta_0^i,\ldots,\theta_n^i$, by setting
$$J_{k,j}^i=\big[ t_{k,j-1},~t_{k,j}\big),\qquad\qquad t_{k,j}=
{T\over N}\cdot\Big( k+\sum_{\ell=0}^j \theta_\ell^i\Big),\qquad
t_{k,-1}={Tk\over N}.$$
For any point $(t,x)\in\overline{\Delta^i}$ we now set
$$\left\{\eqalign{ &g^i(t,x)\doteq\psi^i_j(t,x)\cr
&\bar g^i(t,x)=y_j^i\cr}\right.
\qquad\qquad \hbox{if\ }~t\in \bigcup_{k=1}^N J^i_{k,j}.
\eqno(3.12)$$
The piecewise Lipschitz selection $g$\ and a piecewise constant
approximation $\bar g$ of $g$ can now be defined as
$$g=\sum_{i=1}^{\nu}g^i \chi_{\strut\Delta^i},
\qquad\qquad \bar g=\sum_{i=1}^\nu \bar g^i\chi_{\strut \Delta^i}.
\eqno(3.13)$$
By construction, recalling (3.7) and (3.8),
the conditions (a), (b) in (i) clearly hold.
\v
It remains to show that the estimates in (ii) hold as well.  Let
$\tau,~\tau'\in [0,T]$ and $u\in Y$\ be such that
$\big(t,u(t)\big)\in S$\ for every $t\in [\tau,\tau']$, and define
$$E^i=\left\{t\in I\ ;~~~ (t,u(t))\in\Delta^i\right\},\qquad\qquad i=1,
\ldots,\nu.$$  From our previous definition
$\Delta^i\doteq\Gamma_i\setminus\bigcup_{j<i}\Gamma_j$,
where each $\Gamma_j$ is a $(M+1)$-cone, it follows that every $E^i$\ is the
union of at most $i$\ disjoint intervals.  We can thus write
$$E^i=\Big(\bigcup_{J_k\subset E^i} J_k\Big)\cup \hat{E}^i,$$
with $J_k$ given by (3.11) and
$$m(\hat{E}^i)\leq {2iT\over N}\leq {2\nu T\over N}.\eqno(3.14)$$
Since
$$\phi(t_i,x_i)=\sum_{j=0}^n\theta_j^i\yij,\eqno(3.15)$$
the definition of $\bar g$ at (3.12), (3.13) implies
$$\int_{J_k}\big[\phi\txi-\bar{g}\sus\big]~ds=
m(J_k)\cdot\left[\phi\txi-\sumj \theta_j^i \yij\right]=0.$$
Therefore, from (3.9) and (3.6) it follows
$$\eqalign{&\left|\int_{J_k}\big[\phi\sus-g\sus\big]~ds\right|\leq
\left|\int_{J_k}\big[\phi\sus-\phi\txi\big]~ds\right|\cr
&\qquad +\left|\int_{J_k}\big[\phi\txi-\bar{g}\sus\big]~ds\right|+
     \left|\int_{J_k}\big[\bar{g}\sus-g\sus\big]~ds\right|\cr
 &\quad\leq m(J_k)\cdot\left[{\ve\over 4T}+0+{\ve\over 4T}\right]=
m(J_k)\cdot{\ve\over 2T}.\cr}$$
The choice of $N$ at (3.10) and the bound (3.14) thus imply
$$\left|\int_\tau^{\tau'}\big[\phi\sus-g\sus\big]~ds
\right|\leq
2M\cdot m\Big(\bigcup_{i=1}^\nu\hat{E}^i\Big)+(\tau'-\tau){\ve\over 2T}\leq
2M\nu\cdot{2\nu T\over N}+{\ve\over 2}\leq\ve,$$
proving (3.4).
\v
We next consider (3.5).  For a fixed $i\in\{1,\ldots,\nu\}$, let
$E^i$ be as before and define
$$\xi_{-1}=0,\qquad
\xi_j=\sum_{\ell=0}^j\theta_\ell^i,\qquad\qquad w^i(\xi)=\sumj \yij
\chi_{\strut [\xi_{j-1},~\xi_j]}.$$
Recalling (3.15), the definition of $h$\ at (2.1) and H\"older's inequality
together imply
$$\eqalign{h\big(\phi &\txi,~ F\txi\big)\geq
\left(\int_0^1\big|\phi\txi-w^i(\xi)\big|^2~d\xi\right)\!^{1\over 2}
\cr &\geq \int_0^1\big|\phi\txi-w^i(\xi)\big|~d\xi=
\sumj\theta_j^i\big|\phi\txi-\yij\big|.\cr}$$
Using this inequality we obtain
$$\eqalign{\int_{J_k}\big|\phi\txi-&\bar{g}\sus\big|~ds=
 m(J_k)\cdot \sumj\theta_j^i\big|\phi\txi-\yij\big|\cr
 &\leq m(J_k)\cdot h\big(\phi\txi,F\txi\big)\leq
\eta\cdot m(J_k),\cr}$$
and therefore, by (3.9) and (3.6),
$$\eqalign{\int_{J_k}\big|\phi &\sus-g\sus\big|~ds\cr
& \leq \int_{J_k}\big|\phi\sus-\phi\txi\big|~ds+
\int_{J_k}\big|\bar{g}\sus-g\sus\big|~ds\cr
&\quad +\int_{J_k}\big|\phi\txi-\bar{g}\sus\big|\cr
& \leq m(J_k)\cdot\left[{\ve\over 4T}+{\ve\over 4T}+\eta\right]=
m(J_k)\cdot \left({\ve\over 2T}+\eta\right).\cr}$$
Using again (3.14) and (3.10), we conclude
$$\int_\tau^{\tau'}\big|\phi\sus-g\sus\big|~ds\leq
(\tau'-\tau)\left({\ve\over 2T}+\eta\right)+2M\nu\cdot{2\nu T\over N}\leq
\ve+(\tau'-\tau)\eta.$$\q
\v
Using Lemma 3, given any continuous selection $\tilde f$ of $F$ on $\O$, and
any sequence $(\ve_k)_{k\geq 1}$ of strictly positive numbers, we can
generate a sequence $(f_k)_{k\geq 1}$ of selections from $F$ as follows.

To construct $f_1$, we apply the lemma with $S=\O$, $\phi=f_0$,
$\ve=\ve_1$.  This yields a partition $\big\{ A_1^i;~~i=1,\ldots,
\nu_1\big\}$ of $\O$ and a piecewise Lipschitz selection $f_1$ of $F$
of the form
$$f_1=\sum_{i=1}^{\nu_1} f_1^i\chi_{\strut A_1^i}.$$
\v
In general, at the beginning of the $k$-th step we are given a partition
of $\O$, say $\big\{ A_k^i;~~i=1,\ldots,\nu_k\big\}$, and a
selection
$$f_k=\sum_{i=1}^{\nu_k}f_{k}^i\chi_{\strut A^i_k},$$
where each $f^i_k$ is Lipschitz continuous and satisfies
$$h\big(f_k(t,x),F(t,x)\big)\leq\ve_k\qquad\qquad \forall (t,x)\in
\overline{A_k^i}.$$
We then apply Lemma 3 separately to each $A_k^i$, choosing
$S=\overline{A^i_k}$,~ $\ve=\ve_k$, $\phi=f^i_k$.
This yields a partition $\big\{A^i_{k+1};~~
i=1,\ldots,\nu_{k+1}\big\}$ of $\O$ and functions of the form
$$f_{k+1}=\sum_{i=1}^{\nu_{k+1}} f_{k+1}^i\chi_{\strut A_{k+1}^i},
\qquad\qquad\varphi^i_{k+1}(\cdot)=\langle a_{k+1}^i,\cdot\rangle
+b_{k+1}^i,$$
where each $f_{k+1}^i:\overline{ A_{k+1}^i}\mapsto\rn$ is a Lipschitz
continuous selection from $F$,
satisfying the following estimates:
$$\varphi^i_{k+1}(z)>h\big(z,F(t,x)\big)\qquad\qquad\forall (t,x)\in
A^i_{k+1},\eqno(3.16)$$
$$\varphi^i_{k+1}\big(f^i_{k+1}(t,x)\big)\leq\ve_{k+1}\qquad\qquad
\forall (t,x)\in A^i_{k+1},\eqno(3.17)$$
$$\left|\int_\tau^{\tau'}\big[f_{k+1}\sus-f_k\sus\big]~ds\right|\leq
\ve_{k+1},\eqno(3.18)$$
$$\int_\tau^{\tau'}\big| f_{k+1}\sus-f_k\sus\big|~ds\leq
\ve_{k+1}+\ve_k(\tau'-\tau),\eqno(3.19)$$
for every $u\in Y$ and every $\tau,\tau'$, as long as the values $\sus$
remain inside a single set $A^i_k$, for $s\in [\tau,\tau')$.

Observe that, according to Lemma 3, each $A_k^i$ is closed-open in the
finer topology generated by all $(M+1)$--cones.  Therefore, each
$f_k$ is $\Gamma^{M+1}$--continuous.  By Theorem 2 in [2], the substitution
operator $\S^{f_k}:u(\cdot)\mapsto f_k(\cdot, u(\cdot))$
is continuous from the set
$Y$ defined at (2.5) into ${\cal L}^1\big([0,T];~\rn\big)$.
The Picard map $\P^{f_k}$ is thus continuous as well.

Furthermore, there exists an integer $N_k$ with the following property.
Given any $u\in Y$, there exists a finite partition of $[0,T]$ with nodes
$0=\tau_0<\tau_1<\cdots <\tau_{n(u)}=T$, with $n(u)\leq N_k$,
such that, as $t$ ranges in any $[\tau_{\ell-1},~\tau_\ell)$, the point
$(t,u(t))$ remains inside one single set $A_k^i$.  Otherwise stated, the
number of times in which the curve $t\mapsto (t,u(t))$ crosses
a boundary between two distinct sets $A_k^i,~A_k^j$ is smaller that
$N_k$, for every $u\in Y$.  The construction of the $A_k^i$ in terms of
$(M+1)$--cones implies that all these crossings are transversal.
Since the restriction of $f_k$ to each $A_k^i$ is Lipschitz continuous,
it is clear that every Cauchy problem
$$\dot x(t)=f_k(t,x(t)),\qquad\qquad x(t_0)=x_0$$
has a unique solution,
depending continuously on the initial data $\inda$.

{}From (3.18), (3.19) and the property of $N_k$ it follows
$$\eqalign{ \left|\int_0^t \big[
f_{k+1}\sus-f_k\sus\big]~ds\right|
&\leq\sum_{\ell=1}^{\hat\ell}\left|\int_{\tau_{\ell-1}}^{\tau_\ell}
\big[ f_{k+1}\sus-f_k\sus\big]~ds\right|\cr &\leq
N_k\ve_{k+1},\cr}\eqno(3.20)$$
where $0=\tau_0<\tau_1<\cdots <\tau_{\hat\ell}=t$ are the times at which
the map $s\to\sus$ crosses a boundary between two distinct sets
$A_k^i,~A_k^j$.  Since (3.20) holds for every $t\in [0,T]$, we conclude
$$\big\| \P^{f_{k+1}}-\P^{f_k}\big\|\leq N_k\ve_{k+1}.\eqno(3.21)$$
Similarly, for every $u\in Y$ one has
$$\eqalign{\Big\|f_{k+1}(\cdot,u(\cdot))-f_k(\cdot,u(\cdot))&\Big\|_
{{\cal L}^1([0,T];~\rn)}~ \leq~\sum_{\ell=1}^{n(u)}
\int_{\tau_{\ell-1}}^{\tau_\ell}\big| f_{k+1}\sus-f_k\sus\big|~ds\cr
&\leq\sum_{\ell=1}^{n(u)}\big[\ve_{k+1}+\ve_k(\tau_\ell-\tau_{\ell-1})
\big]~\leq~ N_k\ve_{k+1}+\ve_k T.\cr}\eqno(3.22)$$
Now consider the functions $\varphi_k :\rn\times\O\to\R$, with
$$\varphi_k(y,t,x)\doteq\langle a^i_k,y\rangle +b_k^i\qquad\qquad
\hbox{if}\qquad (t,x)\in A_k^i.\eqno(3.23)$$
{}From (3.16), (3.17) it follows
$$\varphi_k(y,t,x)\geq h(y,F(t,x))\qquad\qquad\forall (t,x)\in \O,~~y\in F(
t,x),\eqno(3.24)$$
$$\varphi_k\big(f_k(t,x),t,x\big)\leq\ve_k\qquad\qquad\forall (t,x)\in\O.
\eqno(3.25)$$
For every $u\in Y$, (3.18) and the linearity of $\varphi_k$
w.r.t. $y$~ imply
$$\eqalign{\bigg|\int_0^T\big[\varphi_k\big( &f_{k+1}\sus,s,u(s)\big)
-\varphi_k\big( f_k\sus,s,u(s)\big)\big]~ds\bigg|\cr
&\leq\sum_{\ell=1}^{n(u)}\max\big\{|a^1_k|,\ldots,|a_k^{\nu_k}|
\big\}\cdot\left|\int_{\tau_{\ell-1}}^{\tau_\ell}
\big[ f_{k+1}\sus-f_k\sus\big]~ds\right|\cr
&\leq N_k\cdot\max\big\{|a_k^1|,\ldots,|a_k^{\nu_k}|\big\}\cdot\ve_{k+1}.
\cr}\eqno(3.26)$$
Moreover, for every $\ell\geq k$, from (3.19) it follows
$$\eqalign{\int_0^T\Big|\varphi_k &\big(f_{\ell+1}\sus,s,u(s)\big)-\varphi_k
\big(f_\ell\sus,s,u(s)\big)\Big|~ds\cr
&\qquad\leq\max\big\{|a^1_k|,\ldots,|a^{\nu_k}_k|\big\}
\cdot\int_0^T\big|f_{\ell+1}\sus-f_\ell\sus\big|~ds\cr
&\qquad\leq\max\big\{|a_k^1|,\ldots,|a_k^{\nu_k}|\big\}\cdot\big(N_\ell
\ve_{\ell+1}+\ve_\ell T\big).\cr}\eqno(3.27)$$
Observe that all of the above estimates hold regardless of the choice of
the $\ve_k$.  We now introduce an inductive procedure for choosing the
constants $\ve_k$, which will yield the convergence of the sequence $f_k$
to a function $f$ with the desired properties.

Given $f_0$ and $\ve_0$, by Lemma 2 there exists $\delta_0>0$ such that, if
$g:\O\mapsto \overline B(0,M)$ and $\big\|\P^g-\P^{f_0}\big\|
\leq\delta_0$, then, for each $\inda$, every solution of (2.7) remains
$\ve_0$-close to the unique solution of (1.3).  We then choose
$\ve_1=\delta_0/2$.

By induction on $k$, assume that the functions $f_1,\ldots, f_k$ have been
constructed, together with the linear functions
$\varphi^i_\ell(\cdot)=\langle a_\ell^i,\cdot\rangle +b_\ell^i$
and the integers $N_\ell$, ~$\ell=1,\ldots,k$.
Let the values $\delta_0,\delta_1,\ldots,\delta_k>0$
be inductively chosen, satisfying
$$\delta_\ell\leq{\delta_{\ell-1}\over 2}\qquad\qquad\ell=1,\ldots,k,
\eqno(3.28)$$
and such that $\big\|\P^g-\P^{f_\ell}\big\|\leq\delta_\ell$
implies that for every $\inda$ the solution set of (2.7) has diameter
$\leq 2^{-\ell}$, for $\ell=1,\ldots,k$.  This is possible again because of
Lemma 2.
For $k\geq 1$ we then choose
$$\ve_{k+1}\doteq \min\left\{{\delta_k\over 2N_k},~{2^{-k}\over N_k},
{}~{2^{-k}\over N_k\cdot\max\big\{|a_\ell^i|;~~1\leq \ell\leq k,~1\leq i\leq
\nu_\ell\big\}}\right\}.\eqno(3.29)$$
\v
Using (3.28), (3.29) in (3.21), with $N_0\doteq 1$, we now obtain
$$\sum_{k=p}^\infty\big\|\P^{f_{k+1}}-\P^{f_k}\big\|
\leq\sum_{k=p}^\infty N_k\cdot{\delta_k\over 2N_k}
\leq\sum_{k=p}^\infty{2^{p-k}\delta_p\over 2}\leq\delta_p
\eqno(3.30)$$
for every $p\geq 0$.  From (3.22) and (3.29) we further obtain
$$\sum_{k=1}^\infty\big\|f_{k+1}(\cdot,u(\cdot))-f_k(\cdot,u(\cdot))
\big\|_{{\cal L}^1}\leq\sum_{k=1}^\infty\left( N_k\cdot{2^{-k}\over N_k}
+{2^{1-k}T\over N_k}\right)\leq\sum_{k=1}^\infty \big(2^{-k}+2^{1-k}T\big)
\leq 1+2T.\eqno(3.31)$$
Define
$$f(t,x)\doteq\lim_{k\to\infty}~f_k(t,x)\eqno(3.32)$$
for all $(t,x)\in\O$ at which the sequence $f_k$ converges.
By (3.31), for every $u\in Y$ the sequence $f_k(\cdot, u(\cdot))$ converges
in ${\cal L}^1\big([0,T];\rn\big)$ and a.e. ~on $[0,T]$.  In particular,
considering the constant functions $u\equiv x\in \overline{B}(D,MT)$, by
Fubini's theorem
we conclude that $f$ is defined a.e. on $\O$. Moreover, the substitution
operators $\S^{f_k}:u(\cdot)\mapsto f_k(\cdot, u(\cdot))$
converge to the operator $\S^f:u(\cdot)\mapsto f(\cdot,u(\cdot))$
uniformly on $Y$.  Since each $\S^{f_k}$ is continuous, $\S^f$
is also continuous. Clearly, the Picard map $\P^f$ is continuous as well.
By (3.30) we have
$$\big\| \P^f-\P^{f_k}\big\|\leq\sum_{k=p}^\infty\big\|\P^{f_{k+1}}-
\P^{f_k}\big\|\leq\delta_p\qquad\qquad\forall p\geq 1.$$
Recalling the property of $\delta_p$, this implies that, for every $p$, the
solution set of (2.7) has diameter $\leq 2^{-p}$.  Since $p$ is arbitrary,
for every $\inda$ the Cauchy problem can have at most one solution.
On the other hand, the existence of such a solution is
guaranteed by Schauder's theorem.  The continuous
dependence of this solution on the initial data $t_0,x_0$,
in the norm of $\AC$, is now
an immediate consequence of uniqueness and of the continuity of the
operators $\S^f,~\P^f$.  Furthermore, for $p=0$, (3.30) yields
$\big\|\P^f-\P^{f_0}\big\|\leq\delta_0.$  The choice of $\delta_0$
thus implies (1.4).
\v
It now remains to prove (1.1).
Since every set $F(t,x)$ is closed, it is clear that $f(t,x)\in F(t,x)$.
For every $u\in Y$ and $k\geq 1$,
by (3.24)--(3.27) the choices of $\ve_k$ at (3.29) yield
$$\eqalign{\int_0^T h\big( &f\sus,~F\sus\big)~ds\leq
\int_0^T\varphi_k\big( f\sus,s,u(s)\big)~ds\cr
&\leq\int_0^T\varphi_k\big(f_k\sus,s,u(s)\big)~ds\cr
&\qquad+\left|\int_0^T \big[\varphi_k\big( f_{k+1}\sus,s,u(s)\big)-
\varphi_k\big(f_k\sus,s,u(s)\big)\big]~ds\right|\cr
&\qquad +\sum_{\ell=k+1}^\infty\int_0^T
\big|\varphi_k\big(f_{\ell+1}(\sus,s,u(s)\big)-\varphi_k\big(
f_\ell\sus,s,u(s)\big)\big|~ds\cr
&\leq 2^{1-k}T+2^{-k}+\sum_{\ell=k+1}^\infty\big( 2^{-\ell}+
2^{1-\ell}T\big).\cr}\eqno(3.33)$$
Observing that the right hand side of (3.33) approaches zero as
$k\to \infty$, we conclude that
$$\int_0^T h\big(f(t,u(t)),~F(t,u(t))\big)~dt=0.$$
By (2.2), given any $u\in Y$, this implies
$f(t,u(t))\in ext F(t,u(t))$ for almost every $t\in [0,T]$.
By possibly redefining $f$ on a set of measure zero,
this yields (1.1). \q
\vsk
\n{\medbf 4 - Applications}
\v
Throughout this section we make the following assumptions.
\v
\item{(H)} $F:\dom\mapsto\overline B(0,M)$\ is a bounded
continuous multifunction
with compact values satisfying (LSP), while $D$ is a compact set such that
$\overline B(D,~MT)\subset\Omega$.
\v
An immediate consequence of Theorem 1 is
\v
\n{\bf Corollary 1.}  {\it Let the hypotheses (H) hold.  Then there exists
a continuous map $(t_0,x_0)\mapsto x(\cdot,t_0,x_0)$ from $[0,T]\times D$
into $\AC$, such that}
$$\left\{\eqalign{ \dot x(t,t_0,x_0)&\in ext F\big(t,x(t,t_0,x_0)\big)
\qquad\qquad\forall t\in [0,T],\cr
x(t_0,t_0,x_0) &=x_0\qquad\qquad\qquad\qquad\forall t_0,x_0.\cr}
\right.$$
\v
Another consequence of Theorem 1 is the contractibility of the sets of
solutions of certain differential inclusions.  We recall here that a metric
space $X$ is contractible if
there exist a point $\tilde u\in X$\ and a continuous mapping
$\Phi:X\times [0,1]\rightarrow X$\ such that:
$$\Phi(v,0)=\tilde{u},\qquad\quad
\Phi(v,1)=v,\qquad\qquad\forall v\in X.$$
The map $\Phi$ is then called a {\sl null homotopy} of $X$.
\v
\n{\bf Corollary 2.}  {\it Let the assumptions (H) hold.
Then, for any $\bar x\in D$, the sets $\M$, $\Me$ of solutions of
$$x(0)=\bar x,\qquad\qquad\dot x(t)\in F(t,x(t))\qquad t\in [0,T],$$
$$x(0)=\bar x,\qquad\qquad \dot x\in ext F(t,x(t))\qquad t\in [0,T],$$
are both contractible in $\AC$.}
\v
\n{\bf Proof.\ } Let $f$\ be a selection from $ext F$~
with the properties stated in Theorem 1.  As usual, we
denote by $x(\cdot,t_0,x_0)$\ the unique solution of the Cauchy problem
(1.2).
Define the null homotopy
$\Phi:{\M}\times [0,1]\rightarrow {\cal M}$\ by setting
$$\Phi(v,\lambda)(t)\doteq\cases{
v(t) &if\qquad $t\in [0,\lambda T]$,\cr
x(t,\lambda T,v(\lambda T))\qquad &if\qquad $ t\in [\lambda T, T]$.
\cr}$$
By Theorem 1, $\Phi$\ is continuous.  Moreover, setting
$\tilde{u}(\cdot)\doteq u(\cdot,0,\bar x)$, we obtain
$$\Phi(v,0)=\tilde{u},\qquad
\Phi(v,1)=v,\qquad\Phi(v,\lambda)\in\M\qquad\forall v\in {\cal M},$$
proving that ${\cal M}$ is contractible.  We now observe that,
if $v\in\Me$, then $\Phi(v,\lambda)\in\Me$ for every $\lambda$.
Therefore, $\Me$ is contractible as well.
\v
Our last application is concerned with feedback controls.
Let $\Omega\subseteq\R^n$ be open, $U\subset\R^m$ compact, and let
$g:[0,T]\times\Omega\times U\rightarrow\rn$\
be a continuous function.
By a well known theorem of Filippov [8],
the solutions of the control system
$$\dot x=g(t,x,u),\qquad\qquad u\in U,\eqno(4.1)$$
correspond to the trajectories of the differential inclusion
$$\dot x\in F(t,x)\doteq\big\{g(t,x,\omega)\;~~~\omega\in U\big\}.
\eqno(4.2)$$
In connection with (4.1), one can consider the ``relaxed" system
$$\dot x=g^{\#}(t,x,u^{\#}),\qquad\qquad u^{\#}\in U^{\#},\eqno(4.3)$$
whose trajectories are precisely those of the differential inclusion
$$\dot x\in F^{\#}(t,x)\doteq \overline{co}F(t,x).$$
The control system (4.3) is obtained defining the compact set
$$U^{\#}\doteq U\times\cdots\times U\times\Delta_n=U^{n+1}\times\Delta_n,$$
where
$$\Delta_n\doteq\left\{\theta=(\theta_0,\ldots,\theta_n)\ ;~~~
\sum_{i=0}^n \theta_i=1,~~ \theta_i\geq 0\quad\forall i\right\}$$
is the standard simplex in $\R^{n+1}$, and setting
$$g^{\#}(t,x,u^{\#})=g^{\#}\big(t,x,(u_0,\ldots,u_n,(\theta_0,\ldots,\theta_n))
\big) \doteq \sum_{i=0}^n \theta_i f(t,x,u_i).$$
Generalized controls of the form $u^{\#}=(u_0,\ldots,u_n,\theta)$\
taking values in the set $U^{n+1}\times\Delta_n$\
are called {\sl chattering controls}.
\v
\n{\bf Corollary 3.}  {\it Consider the control system (4.1),
with $g:[0,T]\times\Omega\times U\mapsto \overline B(0,M)$
Lipschitz continuous.  Let $D$ be a compact set with $\overline B(D;~MT)
\subset\Omega$.
Let $u^{\#}(t,x)\in U^{\#}$\ be a chattering feedback control such that
the mapping
$$(t,x)\mapsto g^{\#}(t,x,u^{\#}(t,x))\doteq f_0(t,x)$$
is Lipschitz continuous.

Then, for every $\varepsilon_0>0$\ there exists a measurable
feedback control $\bar u=\bar u(t,x)$\ with the following properties:
\v
\item{(a)} For every $(t,x)$, one has
$g(t,x,\bar u(t,x))\in extF(t,x)$,~with $F$ as in (4.2).
\v
\item{(b)} for every $\inda$, the Cauchy problem
$$\dot x(t)=g\big(t,x(t),\bar u(t,x(t))\big),\qquad x(t_0)=x_0$$
has a unique solution $x(\cdot,t_0,x_0)$,
\v
\item{(c)} if $y(\cdot,t_0,x_0)$\ denotes the (unique) solution of the
Cauchy problem
$$\dot y=f_0(t,y(t)),\qquad\qquad y(t_0)=x_0,$$
then for every $(t_0,x_0)$\ one has
$$\big| x(t,t_0,x_0)-y(t,t_0,x_0)\big|<\varepsilon_0,\qquad\qquad\forall t
\in [0,T].$$}
\v
\n{\bf Proof.\ }
The Lipschitz continuity of $g$ implies that the multifunction $F$ in (4.2)
is Lipschitz continuous in the Hausdorff metric, hence it satisfies (LSP).
We can thus apply Theorem 1, and obtain a suitable selection
$f$ of $ext F$, in connection with $f_0,~\ve_0$.
For every $(t,x)$, the set
$$W(t,x)\doteq\big\{\omega\in U\ ;~~~ g(t,x,\omega)=f(t,x)\big\}\subset\R^m
$$
is a compact nonempty subset of $U$.
Let $\bar u(t,x)\in W(t,x)$\ be the lexicographic selection.
Then the feedback control $\bar{u}$\ is measurable, and it is
trivial to check that $\bar u$ satisfies all required properties.
\vsk

\centerline{\medbf References}
\v
\i{[1]} J. P. Aubin and A. Cellina, ``Differential Inclusions",
Springer-Verlag, Berlin, 1984.
\v
\i{[2]} A. Bressan, Directionally continuous selections and differential
inclusions, {\it Funkc. Ekvac.} {\bf 31} (1988), 459-470.
\v
\i{[3]} A. Bressan, The most likely path of a differential inclusion,
{\it J. Differential Equations} {\bf 88} (1990), 155-174.
\v
\i{[4]} A. Bressan, Selections of Lipschitz multifunctions generating a
continuous flow, {\it Diff. \& Integ. Equat.} {\bf 4} (1991), 483-490.
\v
\i{[5]} A. Cellina, On the set of solutions to Lipschitzean differential
inclusions, {\it Diff. \& Integ. Equat.} {\bf 1} (1988), 495-500.
\v
\i{[6]} F. S. De Blasi and G. Pianigiani, On the solution set of
nonconvex differential inclusions, preprint.
\v
\i{[7]} A. LeDonne and M. V. Marchi, Representation of Lipschitz
compact convex valued mappings, {\it Rend. Accad. Naz. Lincei}
{\bf 68} (1980), 278-280.
\v
\i{[8]} A. F. Filippov, On certain questions in the theory of optimal
control, {\it SIAM J. Control} {\bf 1} (1962), 76-84.
\v
\i{[9]} A. Ornelas, Parametrization of Caratheodory multifunctions,
{\it Rend. Sem. Mat. Univ. Padova} {\bf 83} (1990), 33-44.
\bye